\documentclass{am2013}
\usepackage{epsfig}
\usepackage{float}
\usepackage{amsmath}
\title{Waves Speed Averaging Impact on \textit{Godunov} type Schemes for Hyperbolic Equations with Discontinuous Coefficients: The linear scalar case}
\usepackage{epstopdf}
\usepackage{graphicx}
\author{L. Remaki 1,2\\
\vskip 2mm {\small
1. BCAM - Basque Centre for Applied Mathematics\\
Mazarredo, 14. 48009 Bilbao Basque Country – Spain\\
2. Department of mathematics ad computer science, Alfaisal University, KSA \\
lremaki@bcamath.org 
}
}

\begin{document}

\maketitle

\begin{abstract}

This paper deals with the waves speed averaging impact impact on \textit{Godunov} type schemes for linear scalar hyperbolic equations with discontinuous coefficients. In many numerical schemes of \textit{Godunov} type used in fluid dynamics, electromagnetic, electro-hydrodynamic problems and so on, usually a \textit{Riemann} problem needs to be solved to estimate fluxes. The exact solution is generally not possible to obtain, but good approximations are provided in many situations like Roe and HLLC \textit{Riemann} solvers in fluids. However all these solvers assume that the acoustic waves speed are continuous by considering some averaging. This could unfortunately lead to a wrong solution as we will show in this paper for the linear scalar case. Providing a \textit{Riemann} solver in the general case of non-linear hyperbolic systems with discontinuous waves speed is a very hard task, therefore in this paper and as a first step, we focus on the linear and scalar case. In a previous work we proposed for such problems a \textit{Riemann} solution that takes into account the discontinuities of the waves speed, we provided a numerical argument to show the validity of the solution. In this paper, first a new argument using regularization technique is provided to reinforce the validity of the proposed solution. Then, the corresponding \textit{Godunov} scheme is derived and the effect of waves speed averaging is clearly demonstrated with a clear connection to the distribution product phenomenon.
\end{abstract}

\section*{keyword}
Hyperbolic equations, \textit{Riemann} solver, waves speed,  \textit{Godunov} scheme, CFD, generalized functions 


\section{Introduction}
 
 In many numerical methods such as finite volume, Discontinuous Galerkin (DG), Discontinuous finite volume, and so on, estimation of numerical fluxes at cells (sub-cells) faces is the most important part of the numerical scheme. The accuracy of the method depends on the accuracy of the flux estimation. For the convective fluxes, generally a \textit{Riemann} problem is considered and then an approximation \textit{Riemann} solver is used. This leads to a stable upwind numerical schemes. This approach was first proposed by \textit{Godunov} \cite{Godunov-1959}, consequently we refer to such methods by \textit{Godunov} type methods. Depending of the problem to solve, many \textit{Riemann} solver approximations where developed. Among the most popular in computational fluid dynamic we can cite the \textit{Roe} solver \cite{Roe-1981, Roe-1986}, the HLL Riemann solver \cite{Harten-1983} and the HLLC solver \cite{Toro-1994}. For the \textit{Roe} solver, the Jacobian matrix is averaged in such a way that hyperbolicity, consistency with the exact Jacobian and conservation across discontinuities still fulfilled. This solver has been modified \cite{Einfeldtoe-1988}, \cite{Einfeldtoe-1991} to overcome the shortcoming for low-density flows. The HLL solver \cite{Harten-1983} solves the original nonlinear flux to take nonlinearity into account. It has a major drawback however due to the space averaging process, the contact discontinuities, shear waves and material interfaces are not captured. To remedy this problem, the HLLC solver \cite{Toro-1994} was proposed by adding the missing wave to the structure. However all these methods assume that the waves speed are continuous across the left and right states of the \textit{Riemann} problem (through the cell interfaces of the mesh) by applying diverse averaging process. This is not true in general; typical situations are recirculation for turbulent flows and transitions from subsonic to supersonic for transonic regimes. The impact of this averaging on the obtained numerical methods has never been addressed. To handle this important issue, we focus in this paper, as a first step, on the scalar and linear case. A \textit{Riemann} solver of scalar hyperbolic linear equation with discontinuous coefficient is developed and published in a proceeding \cite{Rem13} providing a numerical argument of its validity, this was based on a first idea developed in \cite{Rem2}. This solver takes into account the discontinuities of waves speed. To reinforce the validity of the proposed solver and in absence of a rigorous mathematical proof, another strong argument based on a regularization strategy is proposed in this work. Then, and to assess the impact of the averaging process, a \textit{Godunov} type scheme is derived using the proposed \textit{Riemann} solver and compared to the case of averaged waves speed. The test will demonstrate as well that the critical situation is the one where we face a product of distributions, which is not defined in the classical theory of distributions, that could explain the non-validity of the waves speed averaging in such cases.

 \section{\textit{Riemann} Solution for Hyperbolic Equation with Discontinuous Coefficient}
\label{Riemann00}

Consider the following scalar linear hyperbolic equation with discontinuous coefficient,

\begin{equation}\label{eq1}
\begin{array}{rl}
     \frac{\partial}{\partial t}\varphi + a(x)\frac{\partial}{\partial x}\varphi =0, ~~~~\mbox{ on $[o,T]\times \Omega$} \\
     \varphi(0,x) =\varphi_{0} \in L^{\infty}(\Omega)\\

a(.) \in L^{\infty}(\Omega)
 \end{array}
\end{equation}

The initial condition $\varphi_{0}(x)$ and the coefficient $a(x)$ are only bounded functions and can be discontinuous. From the theoretical point of view, the well-posedness of this type of problems is studied in \cite{Rem03}. It is shown that the more critical case is when the solution $\varphi$ and the coefficient $a(.)$ are discontinuous at the same location which leads to a product of distributions (for instance if $a(.)$ is some \textit{Heaviside} function and a Dirac function resulting from the derivative of $\varphi$). This product is not defined in the classical space of distributions which is not an algebra. The well-posdeness of the problem is then studied in a more appropriate space of generalized functions introduced by \textit{J.F Colombeau} , known as well as the \textit{Colombeau} algebra. For more details we refer to \cite{Colom92,ColHei94,Rem03}. We recalled this critical case to make later the connection between this theoretical problem of product of distributions and its impact on the numerical aspect.

Now, let's define the \textit{Riemann} problem associated with problem (\ref{eq1}) 

\begin{equation}\label{eq2}
\begin{array}{rl}
     \frac{\partial}{\partial t}\varphi + a(x)\frac{\partial}{\partial x}\varphi =0, ~~~~\mbox{ on $[o,T]\times \Omega$} \\
     \varphi(0,x) =\varphi_{0}= \left\{ \begin{array}{lr}
     \varphi_{L} &\mbox{ if ~~$x<0$} \\
     \varphi_{R} &\mbox{ if ~~$x>0$}
        \end{array} \right. \\

a(x) = \left\{ \begin{array}{lr}
     a_{L} &\mbox{ if ~~$x<0$} \\
     a_{R} &\mbox{ if ~~$x>0$}
        \end{array} \right.
 \end{array}
\end{equation}
In this equation, the acoustic waves speed $a(.)$ is discontinuous which is not taken into account in the existing \textit{Riemann} solvers where acoustic waves speed are assumed to be continuous in the vicinity of the origin and some averaging is applied to ensure this assumption as we mentioned before. In the proceedings \cite{Rem13}, a \textit{Riemann} solver (a solution to problem (\ref{eq2})) is proposed  based on the following observations of different situations.
\newline
Case 1: $a_{L}>0$ and $a_{R}>0$ we have propagation of the discontinuity (of initial condition) to the right and we do not need to consider what happening within the fan defined by the two acoustic waves, because they will catch up if $a_{L}>\beta_{R}$ and if $a_{L}<\beta_{R}$ an expansion will appear.
\newline
Case 2: $a_{L}<0$ and $a_{R}<0$ similar the previous case with a propagation of the discontinuity to the left.
\newline
Case 3: $a_{L}<0$ and $a_{R}>0$ we have propagation of the discontinuity to the left and the right simultaneously, and we need to determine what happened within the fan defined by the two acoustics waves. We assume that a constant state appears and its expression will be given below.
\newline
Case 4: $a_{L}>0$ and $a_{R}<0$ in this case we have opposite acoustic waves speed and then the discontinuity will remain blocked, which means there is no propagation.

Based on the above observations, the \textit{Riemann} solution of problem (\ref{eq2}) is given by
\begin{equation}\label{eq3}
\varphi (x,t) = \left\{ \begin{array}{rl}
 \varphi_{L} &\mbox{ if ~~$a_{L}>0 ~~and~~ a_{R}>0$} \\
 \lambda &\mbox{ if ~~$a_{L}<0 ~~and~~ a_{R}>0$} \\
\varphi_{R} &\mbox{ if ~~$a_{L}<0 ~~and~~ a_{R}<0$}  \\
 \varphi^{0} &\mbox{ if ~~$a_{L}>0 ~~and~~ a_{R}<0$}
       \end{array} \right.
\end{equation}

Where the expression of the constant $\lambda$ is given by

\begin{equation}\label{eq4}
\begin{array}{rl}

\lambda = \frac{\frac{1}{\vert a_{L} \vert} \varphi_{L} + \frac{1}{\vert a_{R} \vert} \varphi_{R} }{\frac{1}{\vert a_{L} \vert}+ \frac{1}{\vert a_{R} \vert}}

 \end{array}
\end{equation}

\section{On the proof}
In this section we will first recall the numerical argument presented in \cite{Rem13} with more details and then we propose another strong argument that supports the validity of the proposed \textit{Riemann} solver. This argument is based first on a regularization strategy that allows computing the exact solution using characteristics technique and finally obtain the exact solution of the original problem as an asymptotic limit of the regularized solution.

\subsection{Numerical argument}
To prove the validity of solution (\ref{eq3}) and formula (\ref{eq4}) at least numerically, the \textit{Riemann} problem (\ref{eq2}) is solved using a centered second order finite volume scheme stabilized by a first order artificial viscosity, which is equivalent in this case to a finite difference scheme. Let's set \newline
\centerline{ $a^{n}_{i}=a(t_{n},x_{i}) ~~~~~~~ \varphi^{n}_{i}=\varphi(t_{n},x_{i})$ }
Then the explicit numerical scheme is give by 

\centerline{
$
\varphi^{n+1}_{i}= \varphi^{n+1}_{i} + a^{n}_{i}\Delta t(\frac{\varphi^{n}_{i+1}-\varphi^{n}_{i-1}}{2h}) + \varepsilon \Delta t(\frac{\varphi^{n}_{i+1}-2\varphi^{n}_{i}+\varphi^{n}_{i-1}}{h^{2}})
$
}

\begin{figure}[H]
$$\begin{array}{ccc}
    \includegraphics[scale=0.35]{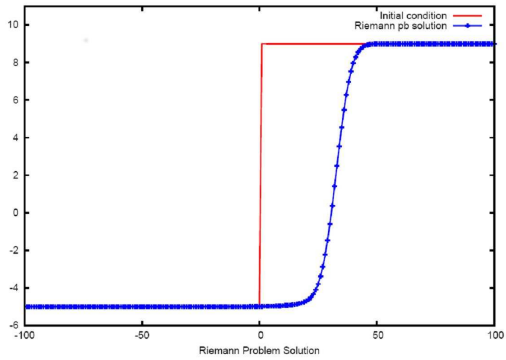} & &\includegraphics[scale=0.35]{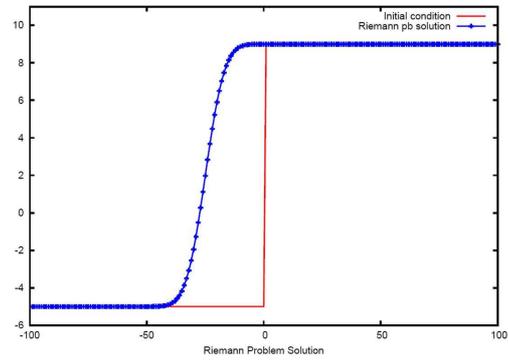} \\
    (a) & & (b) \\
    \includegraphics[scale=0.35]{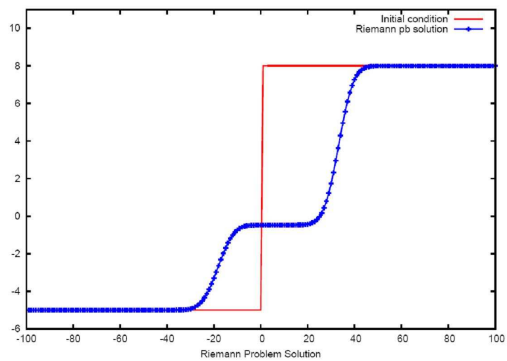} & &\includegraphics[scale=0.35]{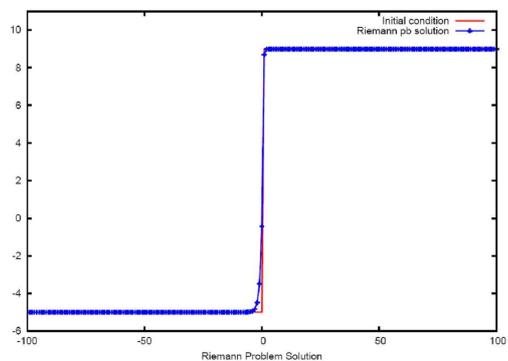} \\
     (c) & & (d)\\
\end{array}$$
\vskip-0.2cm
\caption{ Example 1: Initial condition and \textit{Riemann} solution after 100 time iterations:~~a)~Case (1),~~b)~Case (2),~~c)~Case (3),~~d)~Case (4) }
\label{fig1}
\end{figure}

\begin{figure}[H]
$$\begin{array}{ccc}
    \includegraphics[scale=0.35]{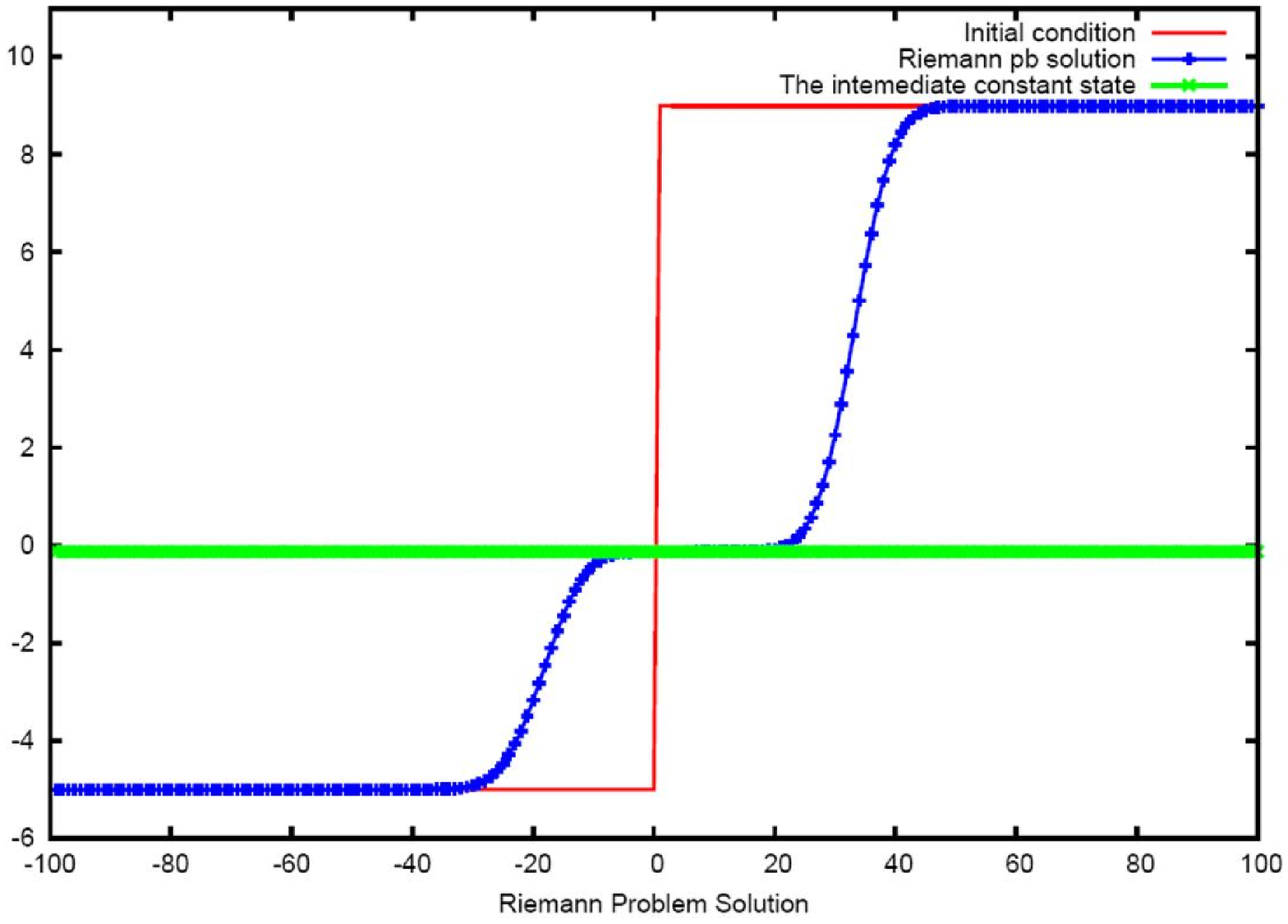} \\
\end{array}$$
\vskip-0.2cm
\caption{ Initial condition and \textit{Riemann} solution corresponding to Case (3) and the value of $\lambda$ given by the proposed formula}
\label{fig2}
\end{figure}

\begin{figure}[h]
$$\begin{array}{ccc}
    \includegraphics[scale=0.35]{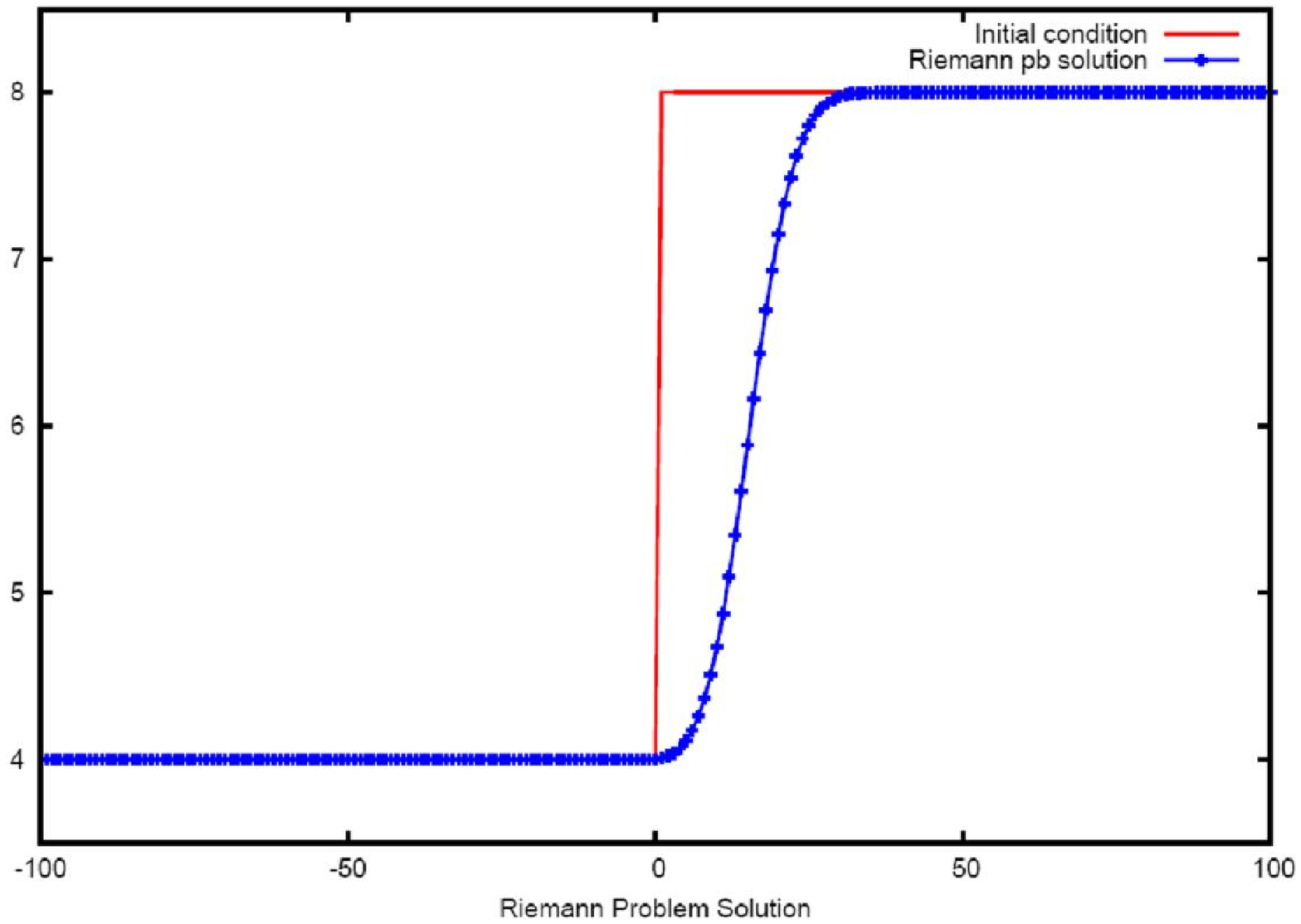} & &\includegraphics[scale=0.35]{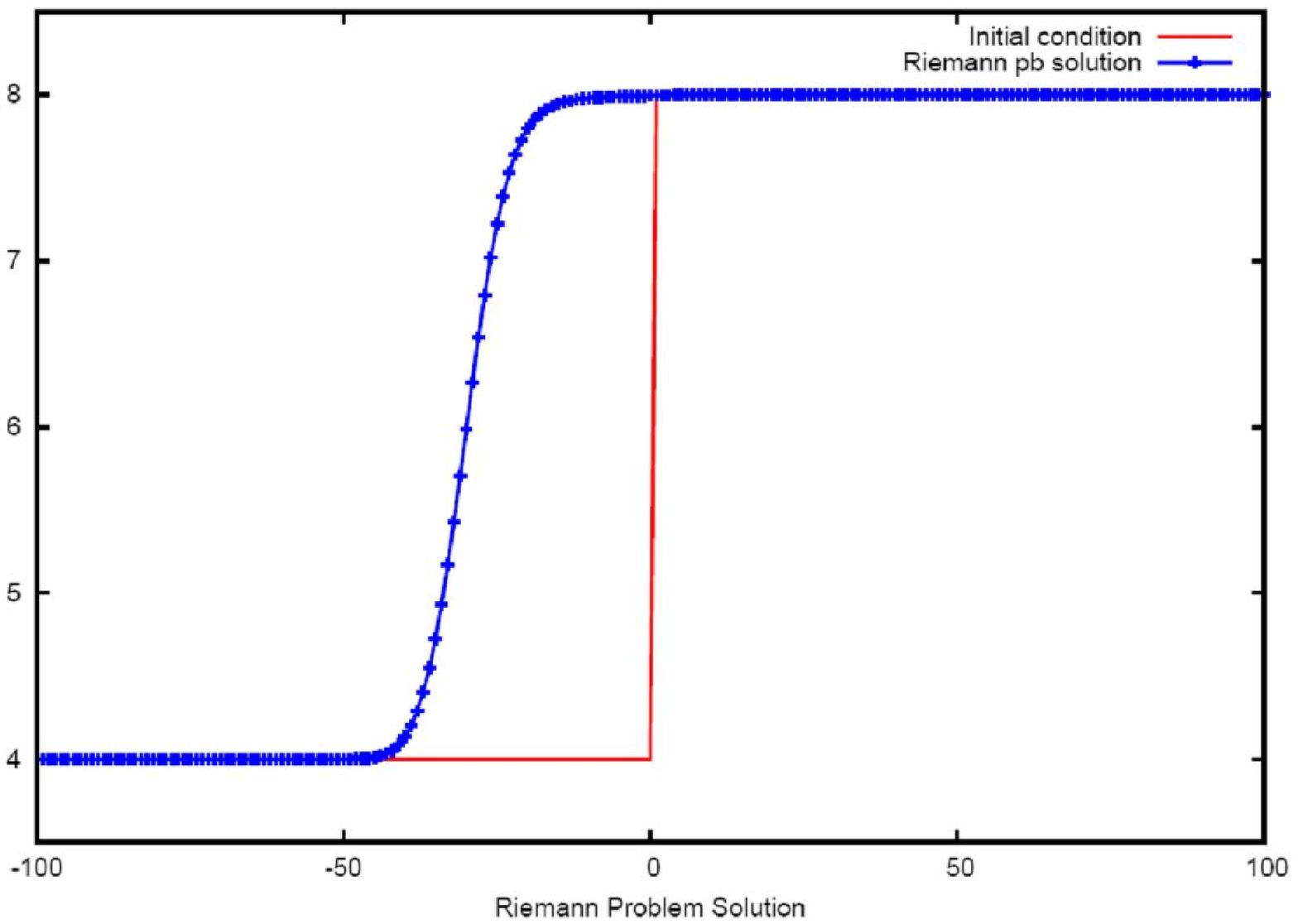} \\
    (a) & & (b) \\
    \includegraphics[scale=0.35]{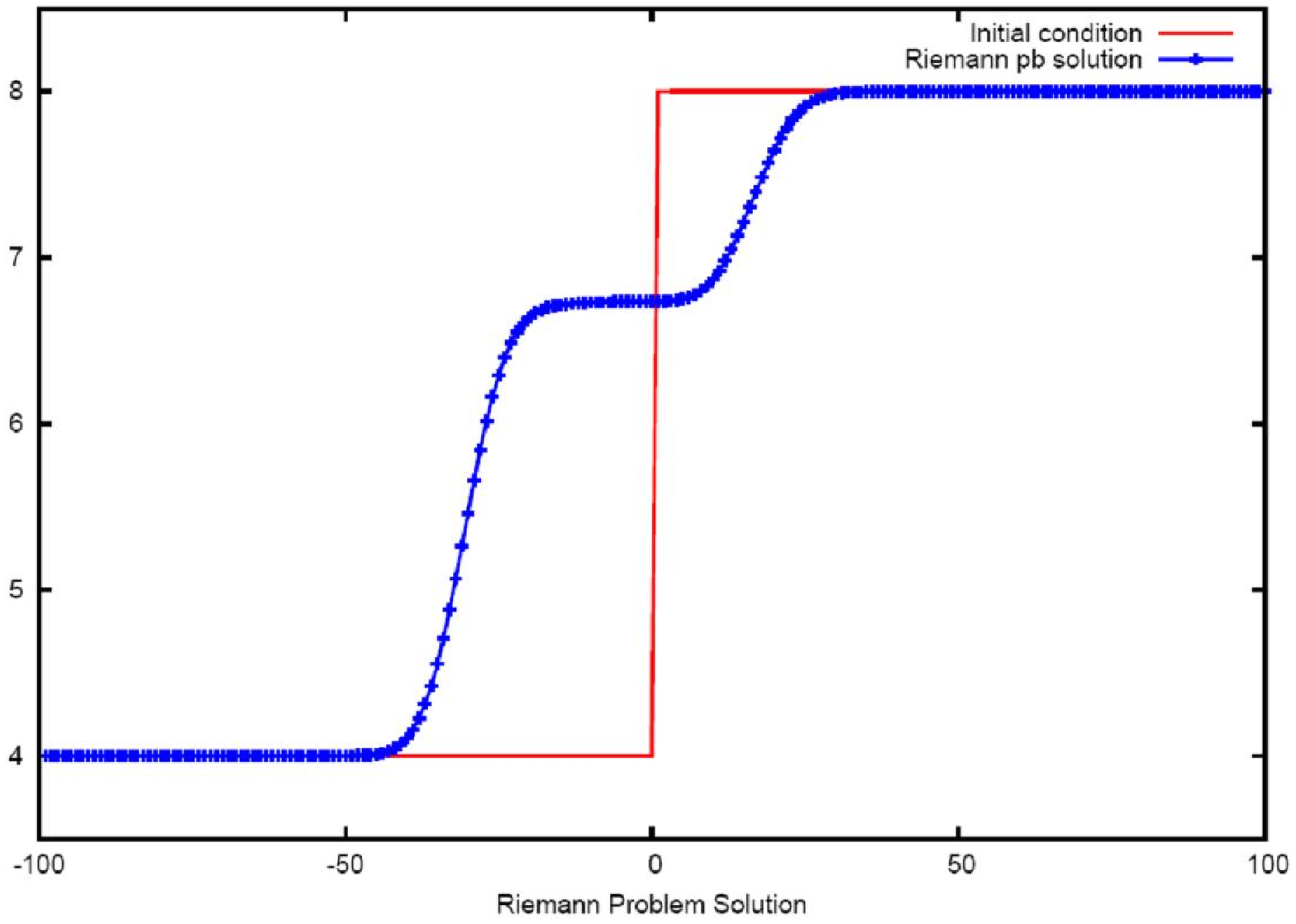} & &\includegraphics[scale=0.35]{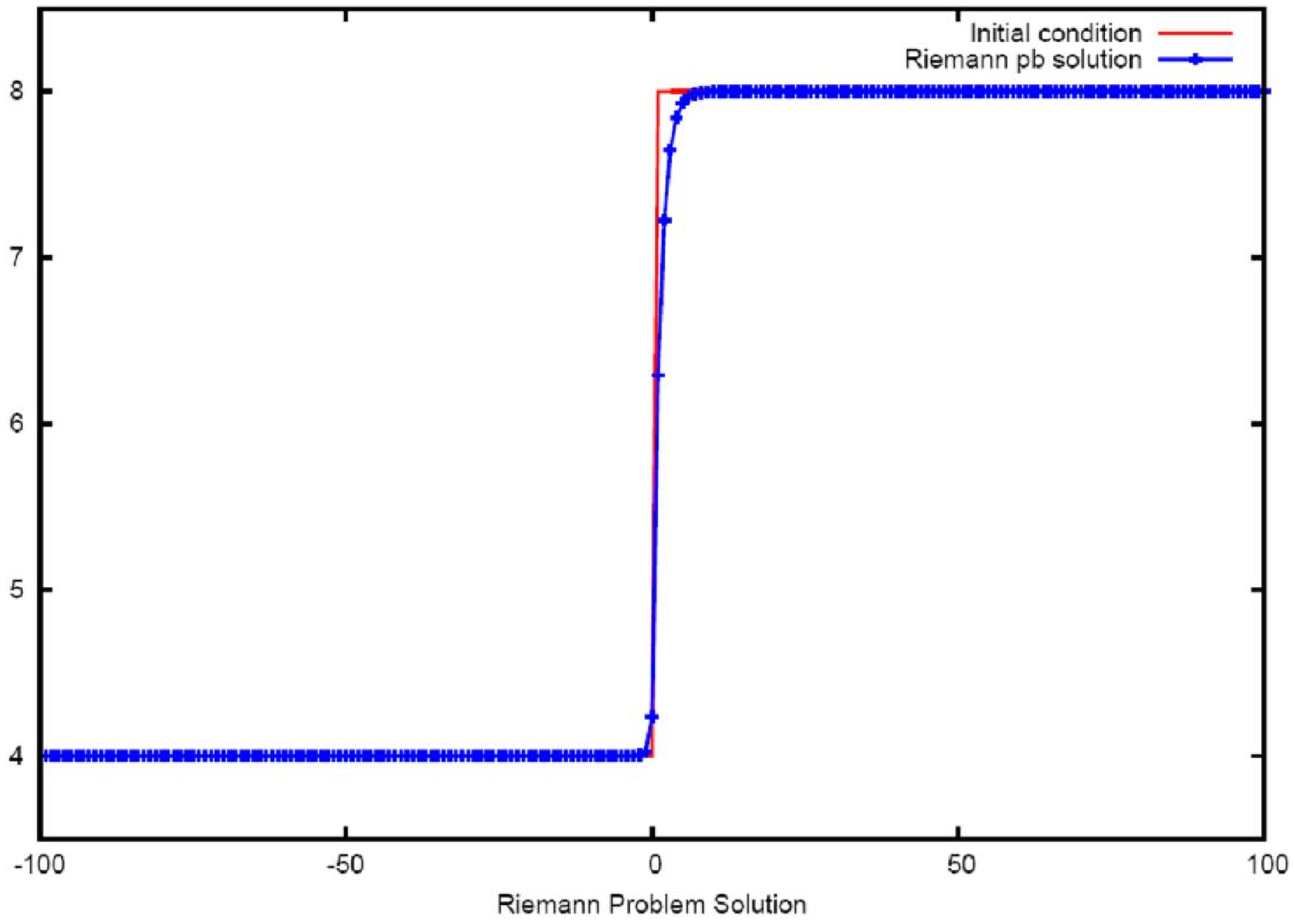} \\
     (c) & & (d)\\
\end{array}$$
\vskip-0.4cm
\caption{ Example 2: Initial condition and \textit{Riemann} solution after 100 time iterations:~~a)~Case (1),~~b)~Case (2),~~c)~Case (3),~~d)~Case (4) }
\label{fig3}
\end{figure}

\begin{figure}[H]
$$\begin{array}{ccc}
    \includegraphics[scale=0.35]{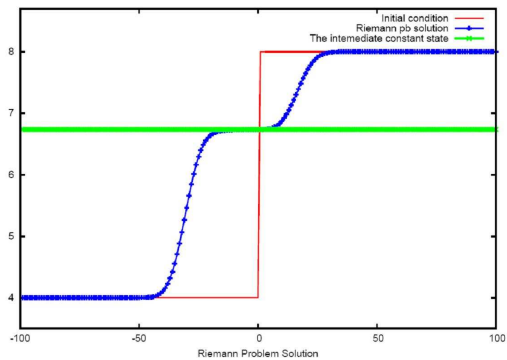} \\
\end{array}$$
\vskip-0.4cm
\caption{Initial condition and \textit{Riemann} solution corresponding to Case (3) and the value of $\lambda$ given by the proposed formula }
\label{fig4}
\end{figure}

Where $\varepsilon$ is the diffusion coefficient that will be set to $0.1$ for the whole simulation.
\vskip+0.2cm
\noindent Several initial conditions $\varphi^{0}$ and acoustic wave speeds $a(.)$ values are tested. All tests demonstrated the validity of the proposed solution. Two examples are shown in Figures (\ref{fig1}) to (\ref{fig4}). We can see in particular the solutions corresponding to cases $3$ and $4$, it is clear that they could not be obtained if the coefficient $a(.)$ is averaged as in the existing \textit{Riemann} solvers. Figures (\ref{fig2}) and (\ref{fig4}) show the perfect agreement of the proposed formula of $\lambda$ with the predicted numerical value.

\subsection{regularization approach argument}
Here we propose to do a regularization of the initial condition $\varphi_{0}$ and the coefficient $a(.)$, in such a way we can compute the exact regularized solution by the characteristics technique, and then the exact solution of the initial problem is obtained by the asymptotic limit.   

Let's replace $\varphi_{0}$ and $a(.)$ by:
\begin{equation}\label{eq5}
\varphi^{\varepsilon}_{0}(x) = \left\{ \begin{array}{ll}
 \varphi_{L} &\mbox{ if ~~$x<-\varepsilon $} \\
 \alpha_{\varepsilon}(x+\varepsilon)^{2}+ \varphi_{L} &\mbox{ if ~~$-\varepsilon \leq x<0 $} \\
-\alpha_{\varepsilon}(x-\varepsilon)^{2}+ \varphi_{R} &\mbox{ if ~~$0\leq x< \varepsilon $}  \\
 \varphi_{R} &\mbox{ if ~~$x \geq \varepsilon$}
       \end{array} \right.
\end{equation}

Where $\alpha_{\varepsilon} = (\varphi_{R} - \varphi_{L})/2\varepsilon^{2}$

\begin{equation}\label{eq6}
a^{\varepsilon}(x) = \left\{ \begin{array}{ll}
 a_{L} &\mbox{ if ~~$x<-\varepsilon $} \\
 \beta_{\varepsilon}(x+\varepsilon)^{2}+ a_{L} &\mbox{ if ~~$-\varepsilon \leq x<0 $} \\
-\beta_{\varepsilon}(x-\varepsilon)^{2}+ a_{R} &\mbox{ if ~~$0 \leq x< \varepsilon $}  \\
 a_{R} &\mbox{ if ~~$x \geq \varepsilon$}
       \end{array} \right.
\end{equation}

Where $\beta_{\varepsilon} = (a_{R} - a_{L})/2\varepsilon^{2}$

The regularization is achieved by connecting the two states by a quadratic functions.
\vskip+0.2cm
Now, consider the same problem with the regularized functions:

\begin{equation}\label{eq7}
\begin{array}{ll}
     \frac{\partial}{\partial t}\varphi^{\varepsilon} + a^{\varepsilon}(x)\frac{\partial}{\partial x}\varphi^{\varepsilon} =0, ~~~~\mbox{ on $[o,T]\times \Omega$} \\
\varphi^{\varepsilon}(0,x) =\varphi^{\varepsilon}_{0}

 \end{array}
\end{equation}

The exact solution is then obtained using the characteristics technique. First we need to solve the equation

\begin{equation}\label{eq8}
\begin{array}{ll}
     y^{'}_{\varepsilon}(s:x,t) = a^{\varepsilon}(s) \\
     y_{\varepsilon}(t) = x
 \end{array}
\end{equation}

Where $y(.)_{\varepsilon}$ is a curve passing by the point $x$ at time $t$

The solution of (\ref{eq7}) is then given by

\begin{equation}\label{eq9}
\varphi^{\varepsilon}(x,t)= \varphi^{\varepsilon}_{0}(y_{\varepsilon}(0:x,t))
\end{equation}

Equation (\ref{eq8}) is solved using the separable variable technique and by distinguishing the intervals $]-\infty, \varepsilon[$, $]-\varepsilon,-\varepsilon + \frac{1}{\sigma_{1}}[$,  $]-\varepsilon + \frac{1}{\sigma_{1}}, 0[$, $]0, \varepsilon - \frac{2}{\sigma_{1}}[$,  $]\varepsilon - \frac{2}{\sigma_{1}}, \varepsilon[$, $]\varepsilon,+\infty[$. The obtained solution is given by

\begin{equation}\label{eq10}
y_{\varepsilon}(s) = \left\{ \begin{array}{ll}
 a_{L}(s-t) +x &\mbox{ if ~~$x<-\varepsilon $} \\
\frac{1-\gamma_{1}e^{2\vartheta_{1}(s-t)}}{\sigma_{1}(1-\gamma_{1}e^{2\vartheta_{1}(s-t)})} -\varepsilon &\mbox{ if ~~$-\varepsilon \leq x< -\varepsilon +\frac{1}{\sigma_{1}} $} \\

\frac{1+ \gamma_{1}e^{2\vartheta_{2}(s-t)}}{\sigma_{1}(1 - \gamma_{1}e^{2\vartheta_{2}(s-t)})} -\varepsilon&\mbox{ if ~~$ -\varepsilon +\frac{1}{\sigma_{1}} \leq x < 0$} \\

\frac{1+ \gamma_{2}e^{2\vartheta_{3}(s-t)}}{\sigma_{1}(-1 + \gamma_{2}e^{2\vartheta_{3}(s-t)})} +\varepsilon&\mbox{ if ~~$ 0 \leq x < \varepsilon -\frac{1}{\sigma_{2}}$} \\

\frac{-1+ \gamma_{2}e^{2\vartheta_{4}(s-t)}}{\sigma_{1}(1 + \gamma_{2}e^{2\vartheta_{4}(s-t)})} +\varepsilon&\mbox{ if ~~$ \varepsilon -\frac{1}{\sigma_{2}} \leq x < \varepsilon$} \\

a_{R}(s-t) +x &\mbox{ if ~~$x\geq \varepsilon $} \\
       \end{array} \right.
\end{equation}

Where $\sigma_{1}=\sqrt{\frac{\beta_{\varepsilon}}{\vert a_{L} \vert}}$, $\sigma_{1}=\sqrt{\frac{\beta_{\varepsilon}}{\vert a_{R} \vert}}$, $\theta_{1}=\sqrt{\vert a_{L}\vert \beta_{\varepsilon}}$, $\theta_{2}=\sqrt{\vert a_{R}\vert \beta_{\varepsilon}}$, $\gamma_{1} = \frac{-\sigma_{1}(\varepsilon+x)+1}{\sigma_{1}(\varepsilon+x)+1}$, $\gamma_{2} = \frac{\sigma_{1}(\varepsilon+x)-1}{\sigma_{1}(\varepsilon+x)+1}$, $\gamma_{3} = \frac{\sigma_{1}(-\varepsilon+x)+1}{\sigma_{1}(-\varepsilon+x)-1}$, $\gamma_{4} = \frac{\sigma_{1}(\varepsilon-x)+1}{\sigma_{1}(-\varepsilon+x)+1}$

Notice that $\gamma_{1}=-\gamma_{2}$ and $\gamma_{3}=-\gamma_{4}$, which implies that expressions of $y_{\varepsilon}(.)$ in the intervals $[-\varepsilon, -\varepsilon +\frac{1}{\sigma_{1}}[$ and $ [-\varepsilon +\frac{1}{\sigma_{1}}, < 0[$ are identical and the same applies to intervals $ [0 ,\varepsilon -\frac{1}{\sigma_{2}}[$ and $[ \varepsilon -\frac{1}{\sigma_{2}} , \varepsilon[$ \\
Therefore the expression of $y_{\varepsilon}(.)$ could be written in a more simple form as:

\begin{equation}\label{eq11}
y(s)_{\varepsilon} = \left\{ \begin{array}{ll}
 a_{L}(s-t) +x &\mbox{ if ~~$x<-\varepsilon $} \\
\frac{1-\gamma_{1}e^{2\vartheta_{1}(s-t)}}{\sigma_{1}(1-\gamma_{1}e^{2\vartheta_{1}(s-t)})} -\varepsilon &\mbox{ if ~~$-\varepsilon \leq x< 0 $} \\

\frac{1+ \gamma_{2}e^{2\vartheta_{3}(s-t)}}{\sigma_{1}(-1 + \gamma_{2}e^{2\vartheta_{3}(s-t)})} +\varepsilon&\mbox{ if ~~$ 0 \leq x < \varepsilon$} \\

a_{R}(s-t) +x &\mbox{ if ~~$x\geq \varepsilon $} \\
       \end{array} \right.
\end{equation}
The exact solution as we mentioned above, is then given by $\varphi^{\varepsilon}(x,t) = \varphi^{\varepsilon}_{0}(y(0:x,t)$.

\begin{figure}[H]
$$\begin{array}{ccc}
    \includegraphics[scale=0.50]{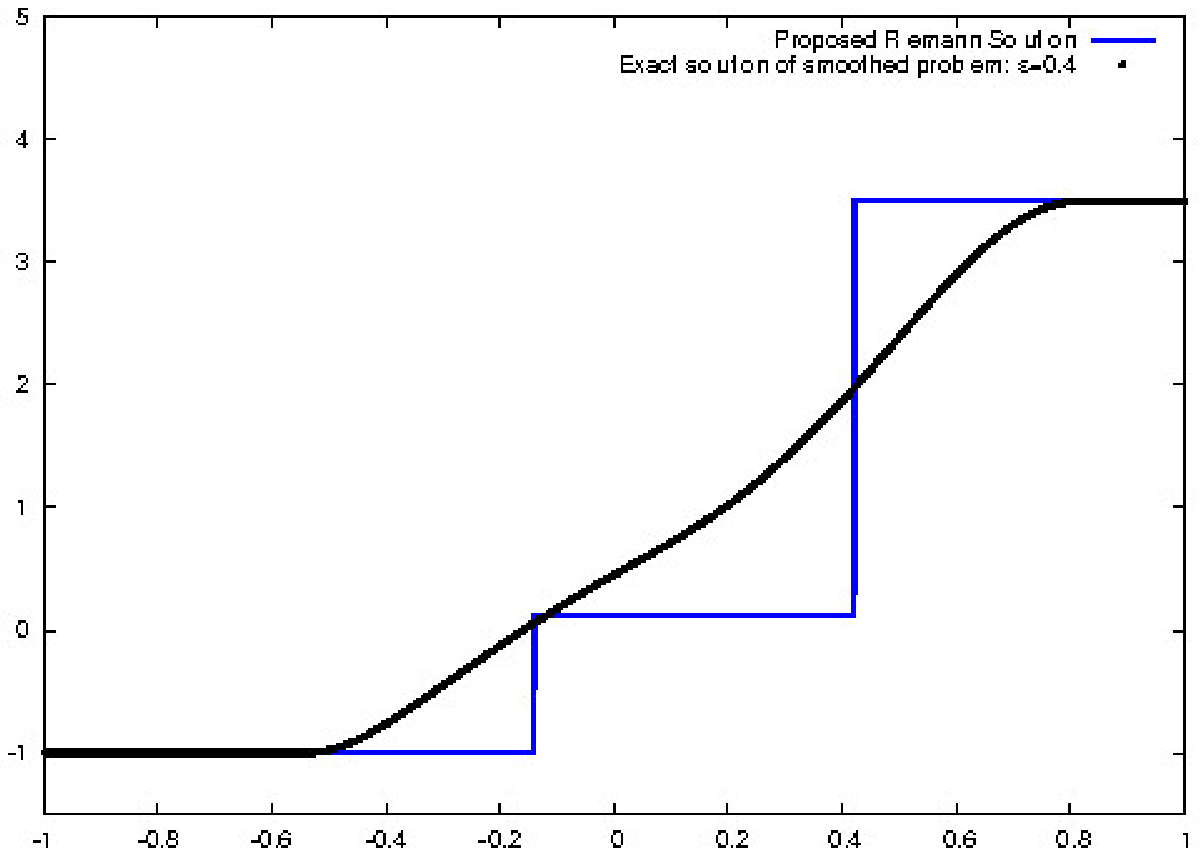} & &\includegraphics[scale=0.50]{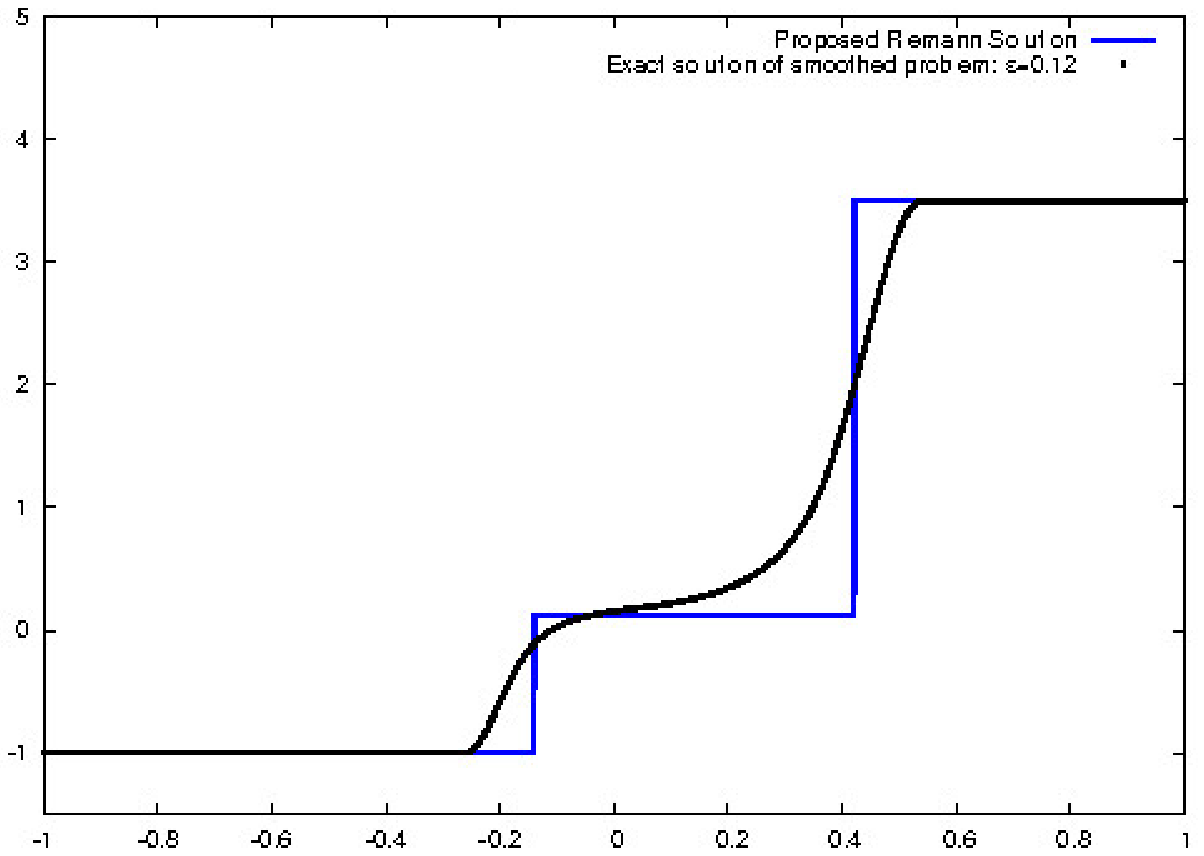}  \\
    (a) & & (b) \\
    \includegraphics[scale=0.50]{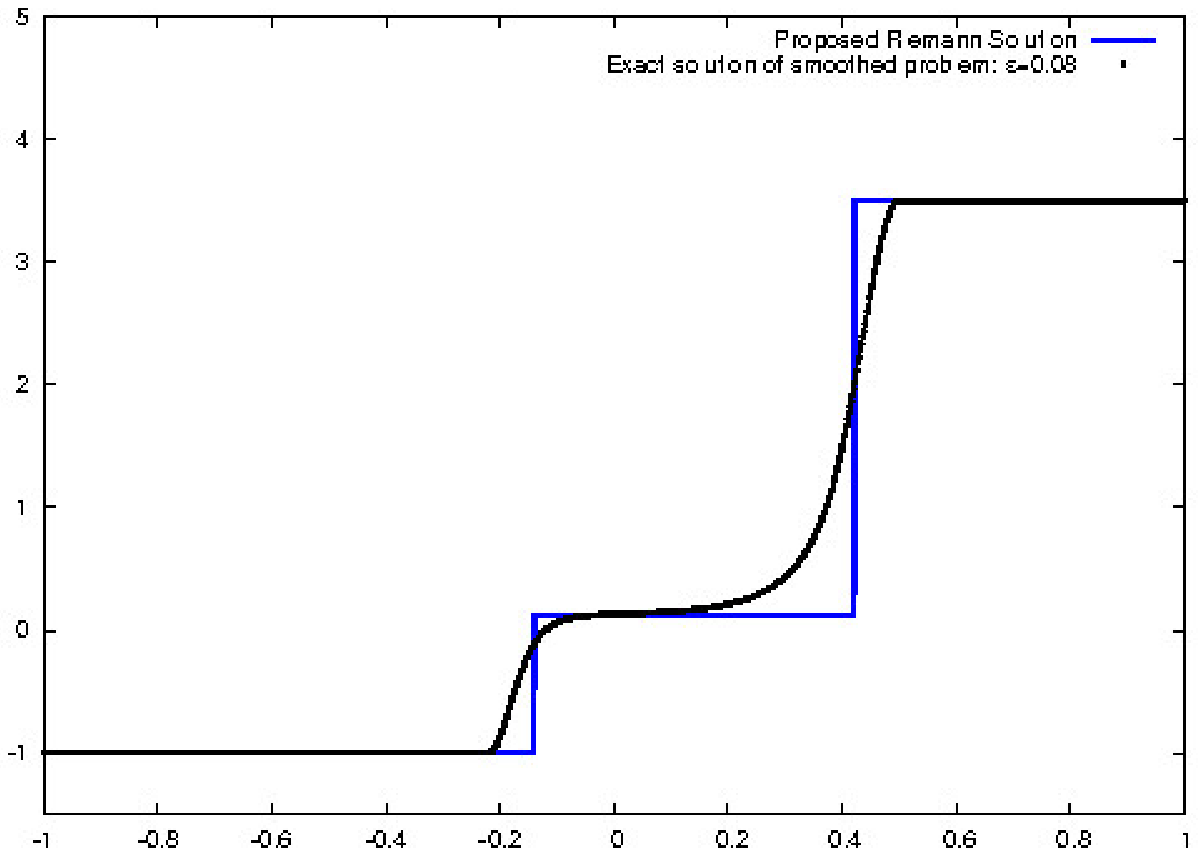} & &\includegraphics[scale=0.50]{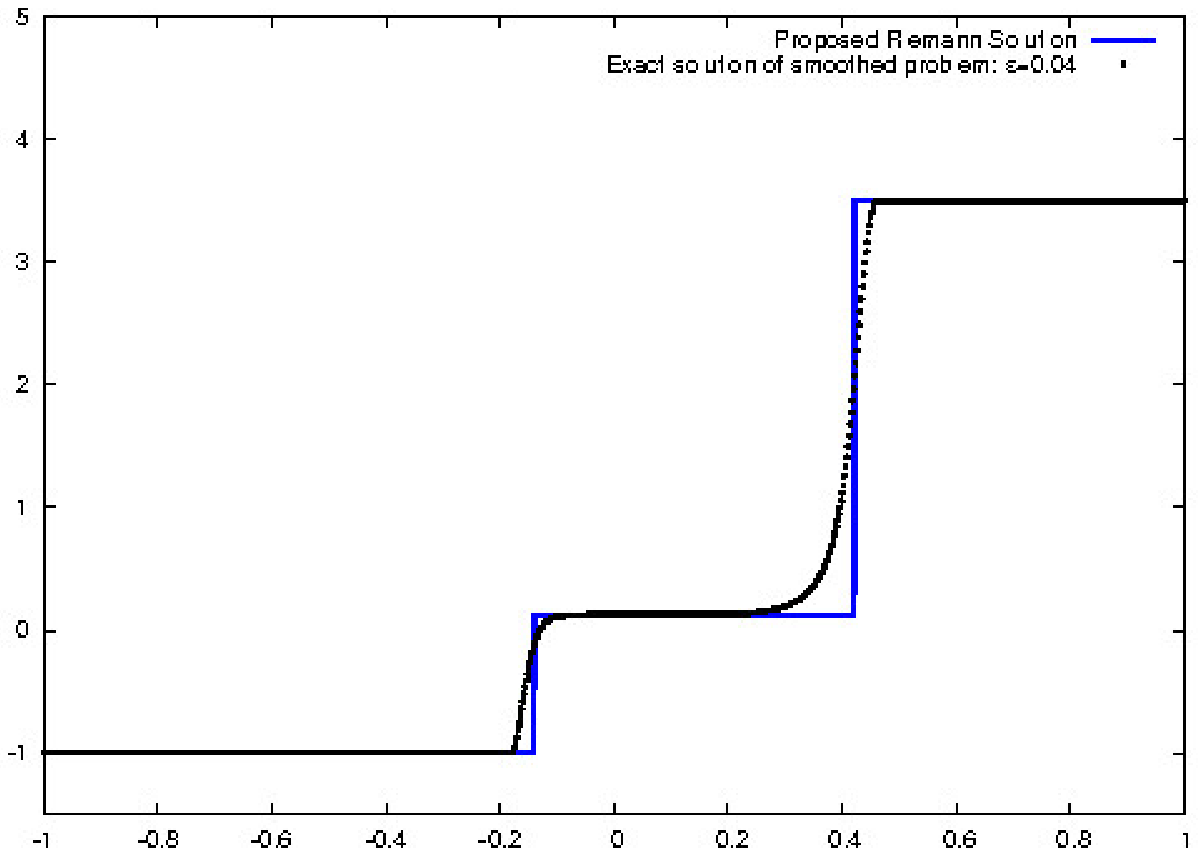}  \\
     (c) & & (d)\\
      \includegraphics[scale=0.50]{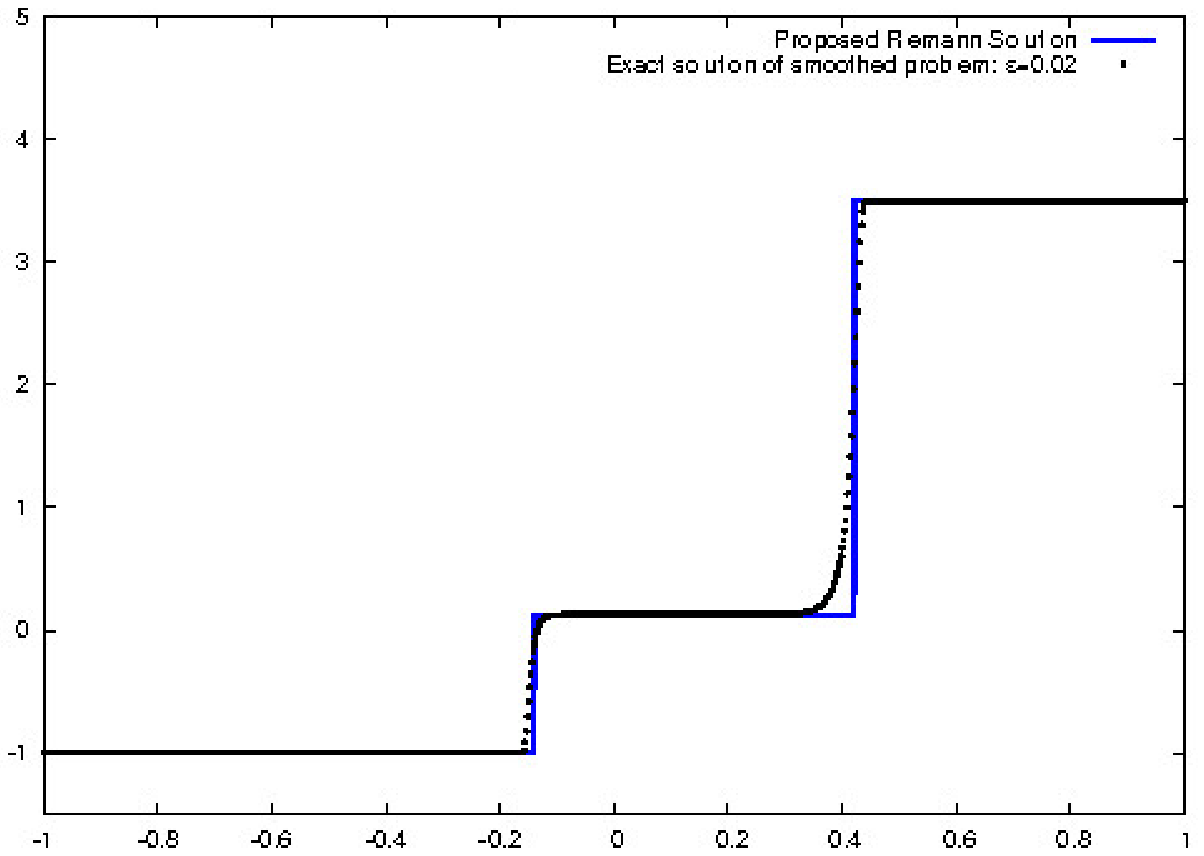} & &\includegraphics[scale=0.50]{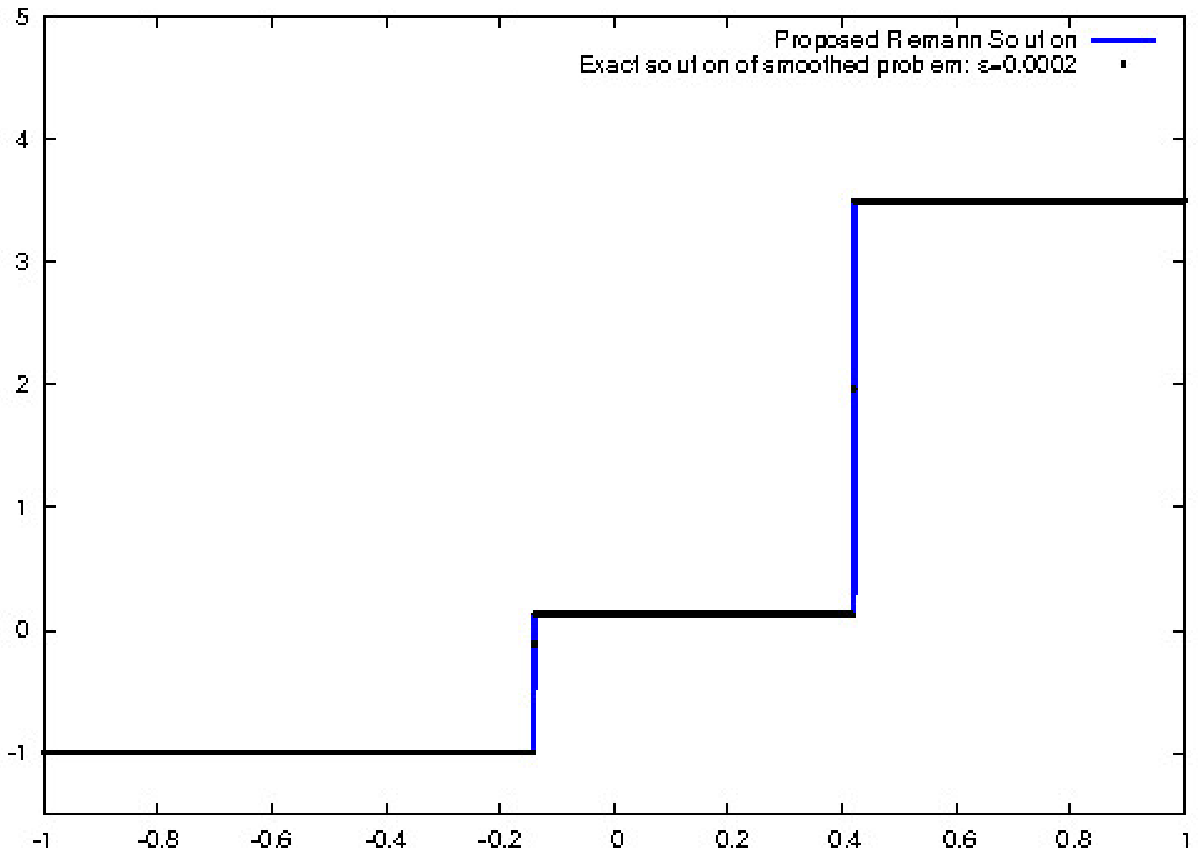}  \\
     (e) & & (f)\\
\end{array}$$
\vskip-0.2cm
\caption{ The proposed \textit{Riemann} solution and the exact regularized solution for different $\varepsilon$ values:~~a)~ $\varepsilon = 0.4$, ~~b)~ $\varepsilon = 0.12$, ~~c)~ $\varepsilon = 0.08$,~~d)~ $\varepsilon = 0.04$,~~e) ~ $\varepsilon = 0.02$,~~f) ~ $\varepsilon = 0.0002$ }
\label{fig5}
\end{figure}

To pass to the limit in $\varepsilon$ to obtain the solution of the initial problem is not easy to perform analytically, this is done numerically instead. The regularized exact solution is implemented and the solution is plotted for different values of $\varepsilon$ and compared to the proposed \textit{Riemann} solver. Note that the implementation is not straightforward, indeed for a given point $(x,t)$ the characteristics passing by this point could jump from and interval to another before crossing the $x$ axes. Therefore we need to pay attention to apply the appropriate formula in (\ref{eq11}). Figure (\ref{fig5}) shows the regularized exact solution for different values of $\varepsilon$. We can see that it is converging to the proposed \textit{Riemann} solver as $\varepsilon$ goes to zero. This demonstrates clearly the validity of the proposed solution.

In the following an answer to the relevant question concerning the impact of waves speed averaging on numerical schemes that use \textit{Riemann} solvers is given. 

\section{Waves speed averaging effect}
In this section we will demonstrate the effect of  waves speed averaging, a process widely used in numericals schemes of hyperbolic equations, by comparing results of a \textit{Godunov} scheme that uses the proposed \textit{Riemann} solver with the one using the averaging process.

\subsection{\textit{Godunov} scheme}
Let's derive the \textit{Godunov} numerical schemes for the linear scalar hyperbolic equations with discontinuous coefficients using the proposed \textit{Riemann} solver and the one for the case of waves speed averaged.  
Let $h$ and $\Delta t$ be the space and time steps, then set $x_{i-1/2}=(i-1/2)h$ and $t_{n}=n\Delta t$, and 

\begin{equation}\label{eq12}
a^{n}_{i}=a(t_{n},x_{i}) ~~~~~~~
\varphi^{n}_{i}=\varphi(t_{n},x_{i})
\end{equation}

The \textit{Godunov} scheme approximation is then given by 
\begin{equation}\label{eq13}
\varphi^{n+1}_{i}=(\varphi^{n+1,R}_{i-1/2}+\varphi^{n+1,L}_{i+1/2})/2
\end{equation}
where $\varphi^{n+1,L}_{i-1/2}$ and $\varphi^{n+1,R}_{i-1/2}$ are obtained by the classical projection process and are given in the following for both sachems.

\subsubsection{Non-averaging waves speed case}
The terms above are given by the following formulas using the propose \textit{Riemann} solver and depending on the sign of $a^{n}_{i-1}$ and $a^{n}_{i}$
\vskip+0.2cm
Case $a^{n}_{i-1} <0$ and $a^{n}_{i} >0$
\begin{figure}[H]
$$\begin{array}{ccc}
    \includegraphics[scale=0.35]{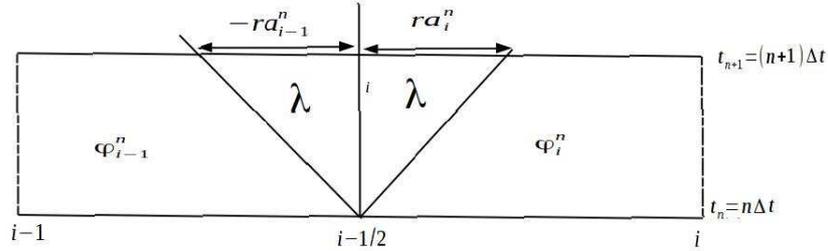} \\
\end{array}$$
\vskip-0.4cm
\caption{\textit{Riemann} Solution:~Case of intermediate constant }
\label{fig6}
\end{figure}

\begin{equation}\label{eq14}
F_{1} = \left\{ \begin{array}{ll}
\varphi^{n+1,R}_{i-1/2}=\varphi^{n}_{i} -2ra^{n}_{i}(-\lambda^{n}_{i-1/2} + \varphi^{n}_{i}) \\
\varphi^{n+1,L}_{i-1/2}=\varphi^{n}_{i-1} -2ra^{n}_{i-1}(\lambda^{n}_{i-1/2} + \varphi^{n}_{i-1})
\end{array} \right.
\end{equation}

with
\begin{equation}\label{eq15}
\lambda^{n}_{i-1/2} = \frac{\frac{1}{\vert a^{n}_{i-1}\vert} \varphi^{n}_{i-1} + \frac{1}{\vert a^{n}_{i}\vert} \varphi^{n}_{i} }{\frac{1}{\vert a^{n}_{i-1}\vert}+ \frac{1}{\vert a^{n}_{i}\vert}}\nonumber
\end{equation}

Case $a^{n}_{i-1} <0$ and $a^{n}_{i} <0$
\begin{figure}[H]
$$\begin{array}{ccc}
    \includegraphics[scale=0.35]{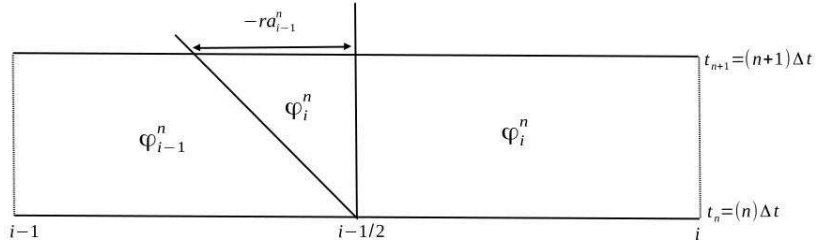} \\
\end{array}$$
\vskip-0.4cm
\caption{\textit{Riemann} Solution:~Case of negative speed }
\label{fig7}
\end{figure}
\begin{equation}\label{eq16}
F_{2} = \left\{\begin{array}{ll}
\varphi^{n+1,R}_{i-1/2}=\varphi^{n}_{i} \\
\varphi^{n+1,L}_{i-1/2}=\varphi^{n}_{i-1} -2ra^{n}_{i-1}(\varphi^{n}_{i} - \varphi^{n}_{i-1})
\end{array} \right.
\end{equation}

Case $a^{n}_{i-1} >0$ and $a^{n}_{i} <0$
\begin{figure}[H]
$$\begin{array}{ccc}
    \includegraphics[scale=0.35]{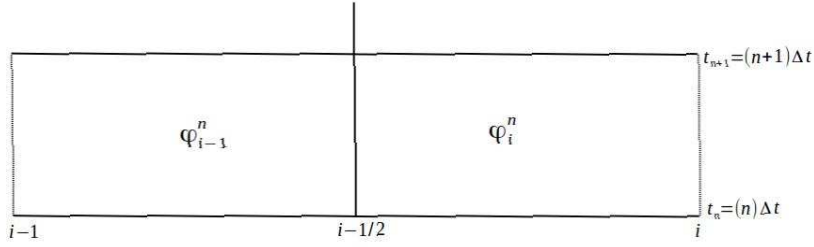} \\
\end{array}$$
\vskip-0.4cm
\caption{\textit{Riemann} Solution:~Case of blocked wave}
\label{fig8}
\end{figure}

\begin{equation}\label{eq17}
F_{3} = \left\{\begin{array}{ll}
\varphi^{n+1,R}_{i-1/2}=\varphi^{n}_{i} \\
\varphi^{n+1,L}_{i-1/2}=\varphi^{n}_{i-1}
\end{array} \right.
\end{equation}

Case $a^{n}_{i-1} >0$ and $a^{n}_{i} >0$
\begin{figure}[H]
$$\begin{array}{ccc}
    \includegraphics[scale=0.35]{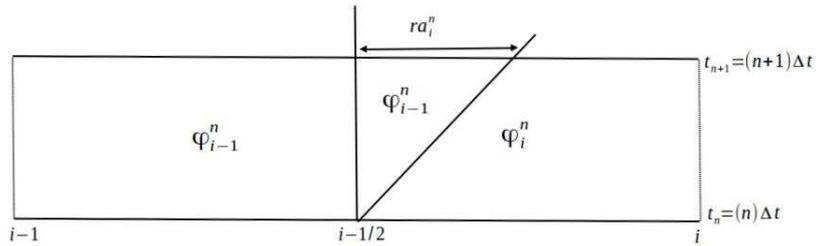} \\
\end{array}$$
\vskip-0.4cm
\caption{\textit{Riemann} Solution:~Case of positive speed }
\label{fig9}
\end{figure}
\begin{equation}\label{eq18}
F_{4} = \left\{ \begin{array}{ll} \varphi^{n+1,R}_{i-1/2}=\varphi^{n}_{i} -2ra^{n}_{i}(\varphi^{n}_{i} - \varphi^{n}_{i-1}) \\
\varphi^{n+1,L}_{i-1/2}=\varphi^{n}_{i-1}
\end{array} \right.
\end{equation}

\vskip+0.2cm
\subsubsection{Averaging waves speed case}
\noindent
For the case where the waves speed are averaged, the \textit{Godunov} scheme is obtained in the same way with formulas reduced to only two cases as a result of the averaging process, namely we have;

\begin{equation}\label{eq19}
\begin{split}
Case ~~ of~~ \alpha^{n}_{i-\frac{1}{2}} = 0.5(a^{n}_{i}+a^{n}_{i-1})>0   \\
F_{1} = \left\{ \begin{array}{ll} \varphi^{n+1,R}_{i-1/2}=\varphi^{n}_{i} -2r\alpha^{n}_{i-\frac{1}{2}}(\varphi^{n}_{i} - \varphi^{n}_{i-1}) \\
\varphi^{n+1,L}_{i-1/2}=\varphi^{n}_{i-1}
\end{array} \right.
\end{split}
\end{equation}

\begin{equation}\label{eq20}
\begin{split}
Case ~~ of~~ \alpha^{n}_{i-\frac{1}{2}} = 0.5(a^{n}_{i}+a^{n}_{i-1})<0 \\
F_{2} = \left\{\begin{array}{ll}
\varphi^{n+1,R}_{i-1/2}=\varphi^{n}_{i} \\
\varphi^{n+1,L}_{i-1/2}=\varphi^{n}_{i-1} -2r\alpha^{n}_{i-\frac{1}{2}}(\varphi^{n}_{i} - \varphi^{n}_{i-1})
\end{array} \right.
\end{split}
\end{equation}

\subsection{Numerical tests}

\begin{figure}[H]
$$\begin{array}{ccc}
    \includegraphics[scale=0.55]{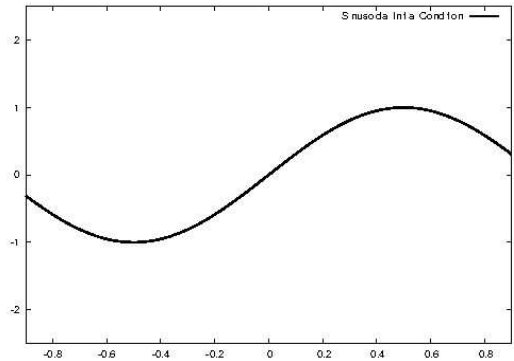} \\
\end{array}$$
\vskip-0.4cm
\caption{Sinusoidal initial condition }
\label{fig10}
\end{figure}

\begin{figure}[H]
$$\begin{array}{ccc}
    \includegraphics[scale=0.55]{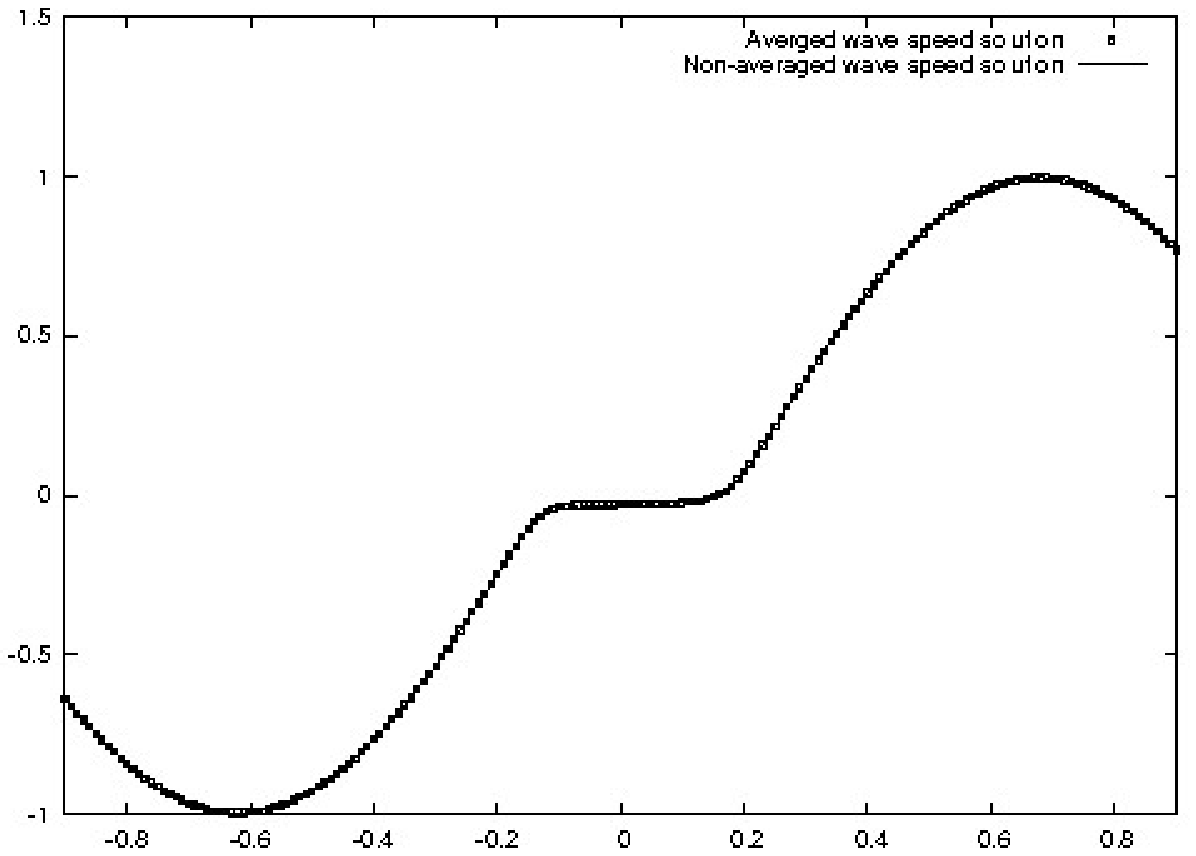} && \includegraphics[scale=0.55]{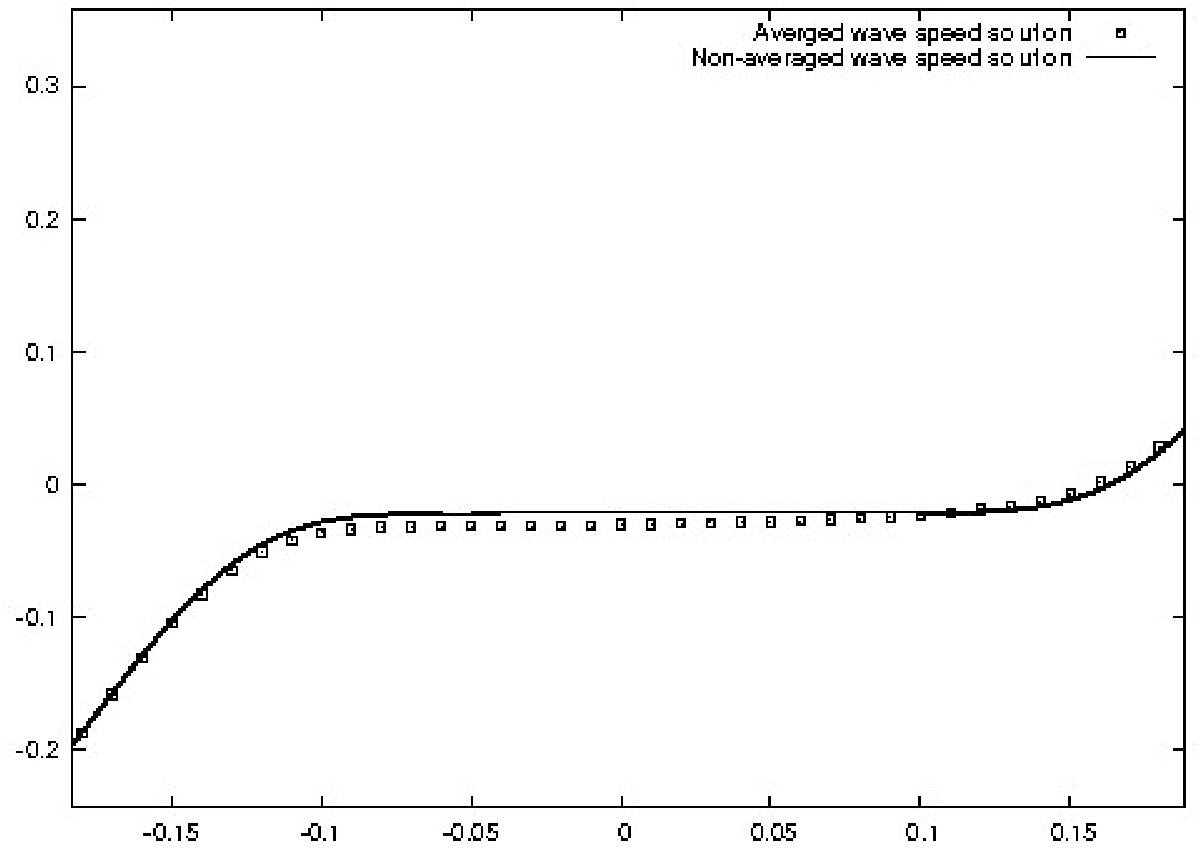} \\
     (a) & & (b) \\
\end{array}$$
\vskip-0.4cm
\caption{(a) Propagation of the initial sinusoidal condition using the proposed \textit{Riemann} solver and the averaged on, ~(b) Zoom on the plateau (the intermediate constant)}
\label{fig11}
\end{figure}

To demonstrate the effect of waves speed averaging we first consider the coefficient $a(.)$ to be a two states function and set to 
\begin{equation}
a(x) =\left\{\begin{array}{ll}
-2      ~~~~~   if x < 0 \\
~~3       ~~~~~  if x > 0
\end{array} \right.
\end{equation}
Which means that we only considering the case of of the appearing the of intermediate constant state $\lambda$ in the proposed \textit{Riemann} solver. As for the initial condition, first we consider a sinusoidal function, in this case we avoid to have a product of distribution of type \textit{Heaviside} times a Dirac function (not defined in the classical space of distributions) which was the main motivation to the use of the \textit{Colombeau} algebra for the mathematical analyse of the problem. Results of the test are depicted in Figures (\ref{fig10},\ref{fig11}). As we can see both schemes give the same results, this means that in this case averaging waves speed allows to well capture the intermediate state with a correct value.

\begin{figure}[H]
$$\begin{array}{ccc}
    \includegraphics[scale=0.55]{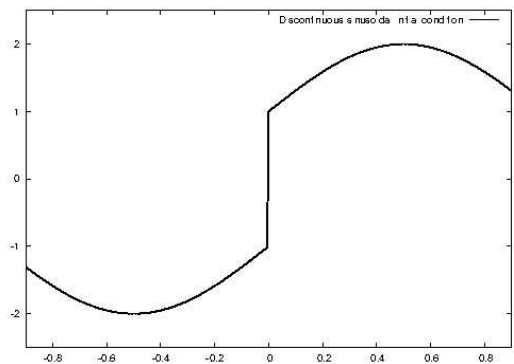} \\
\end{array}$$
\vskip-0.4cm
\caption{Discontinuous sinusoidal initial condition }
\label{fig12}
\end{figure}

\begin{figure}[H]
$$\begin{array}{ccc}
    \includegraphics[scale=0.55]{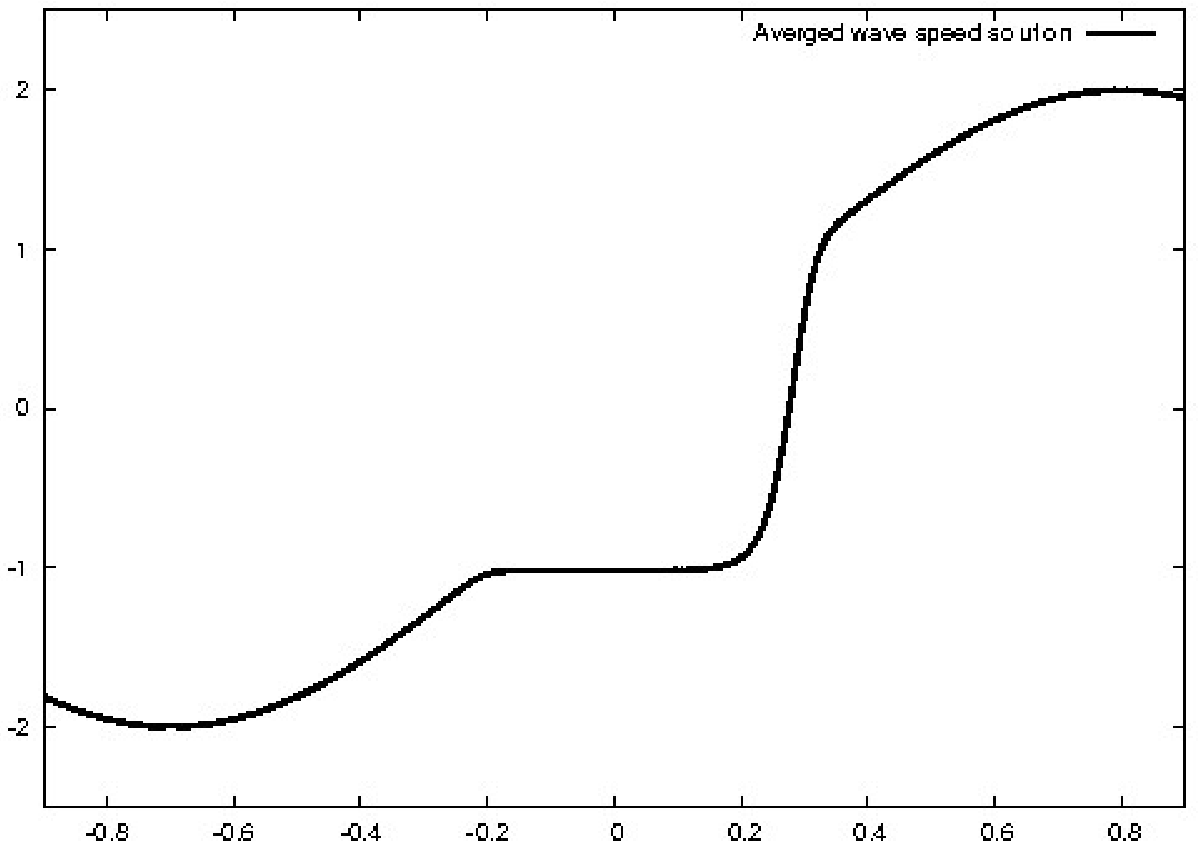} && \includegraphics[scale=0.55]{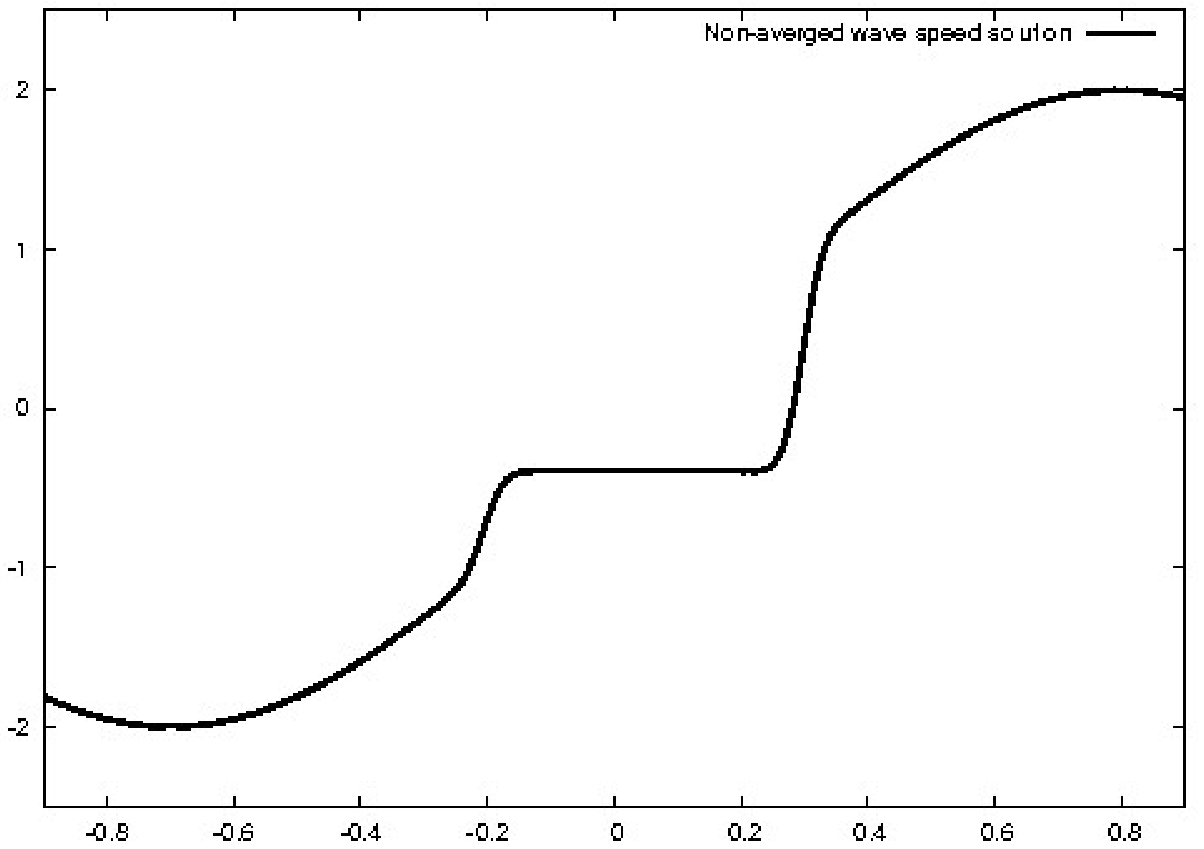} \\
     (a) & & (b) \\
\end{array}$$
\vskip-0.4cm
\caption{(a) Propagation of the initial sinusoidal condition using the averaged \textit{Riemann} solver, ~(b) Propagation of the initial sinusoidal condition using the proposed \textit{Riemann} solver }
\label{fig13}
\end{figure}

\begin{figure}[H]
$$\begin{array}{ccc}
    \includegraphics[scale=0.55]{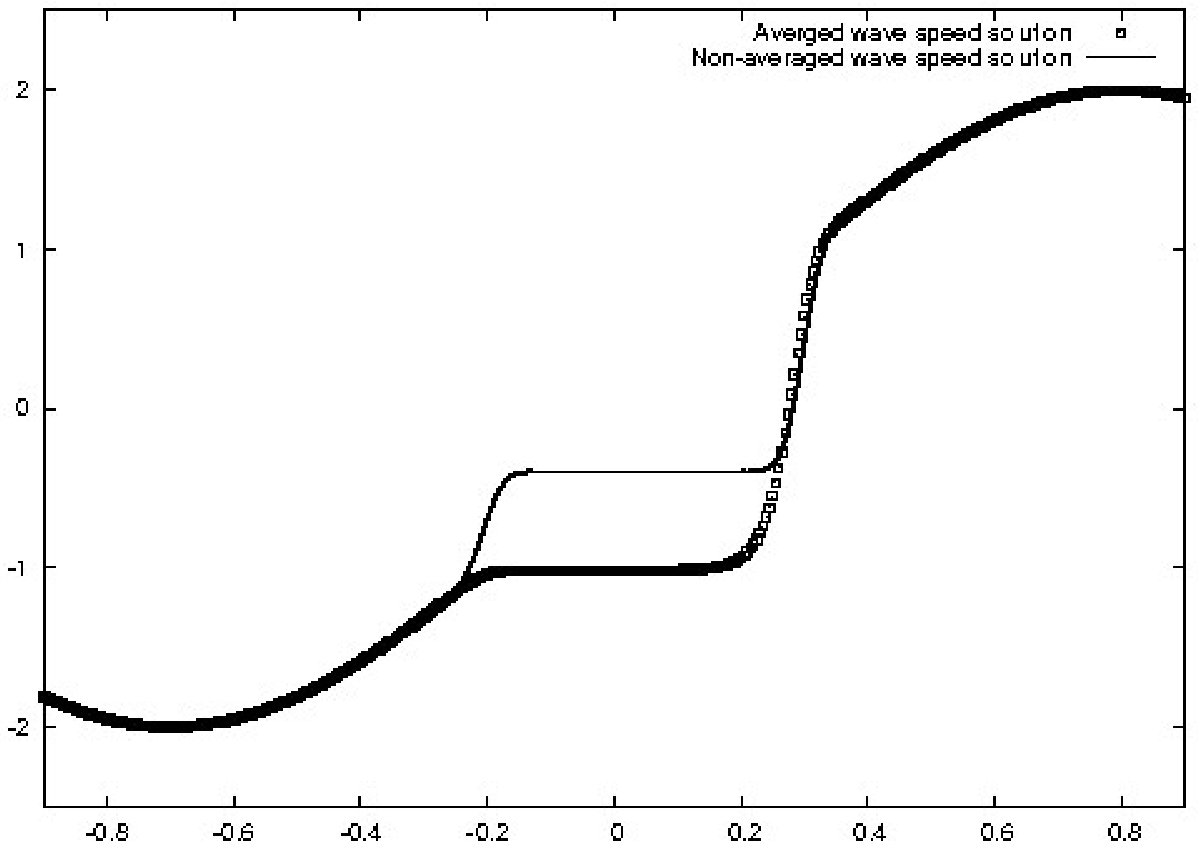} \\
\end{array}$$
\vskip-0.4cm
\caption{Curves of Figure (\ref{fig13}) superimposed}
\label{fig14}
\end{figure}

\begin{figure}[H]
$$\begin{array}{ccc}
    \includegraphics[scale=0.55]{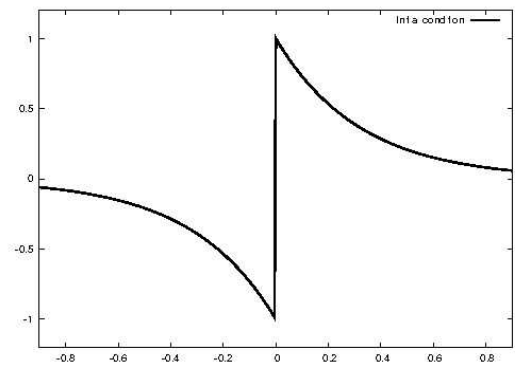} \\
\end{array}$$
\vskip-0.4cm
\caption{Discontinuous polynomial initial condition}
\label{fig15}
\end{figure}

\begin{figure}[H]
$$\begin{array}{ccc}
    \includegraphics[scale=0.55]{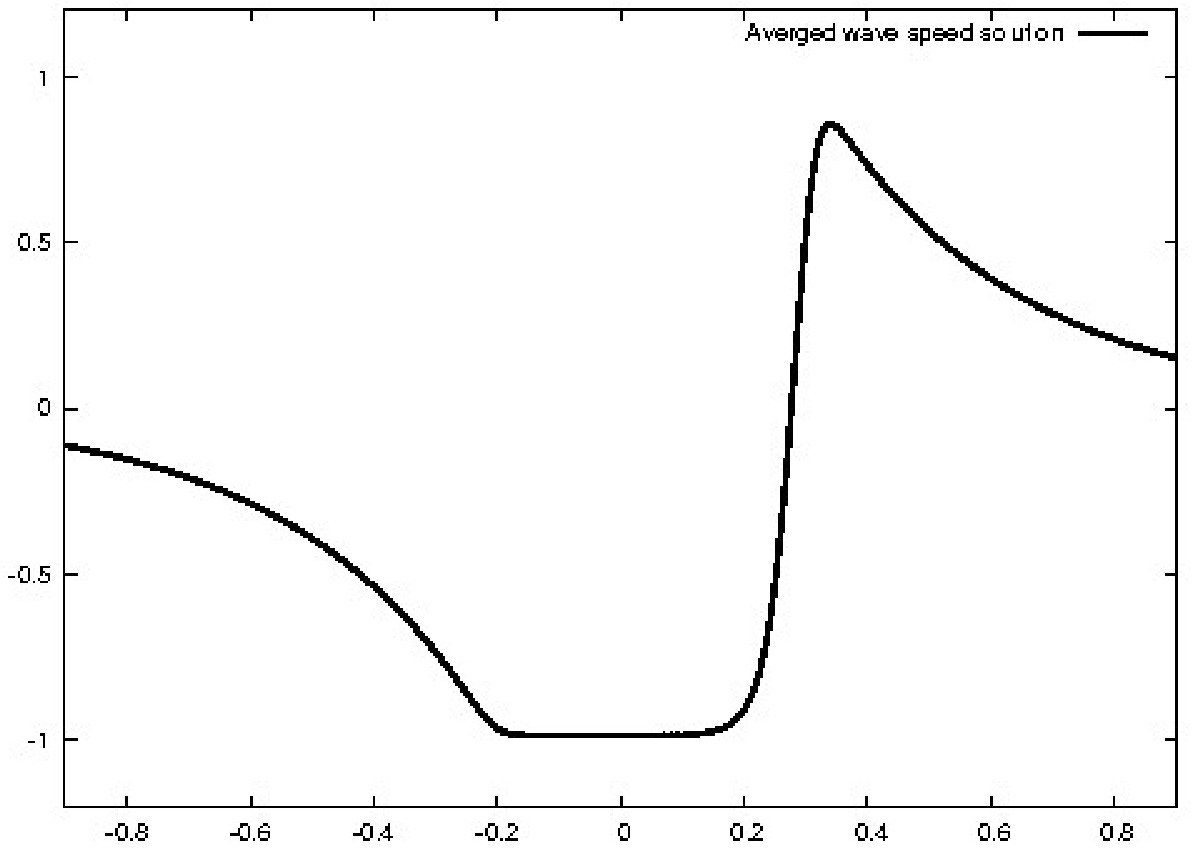} && \includegraphics[scale=0.55]{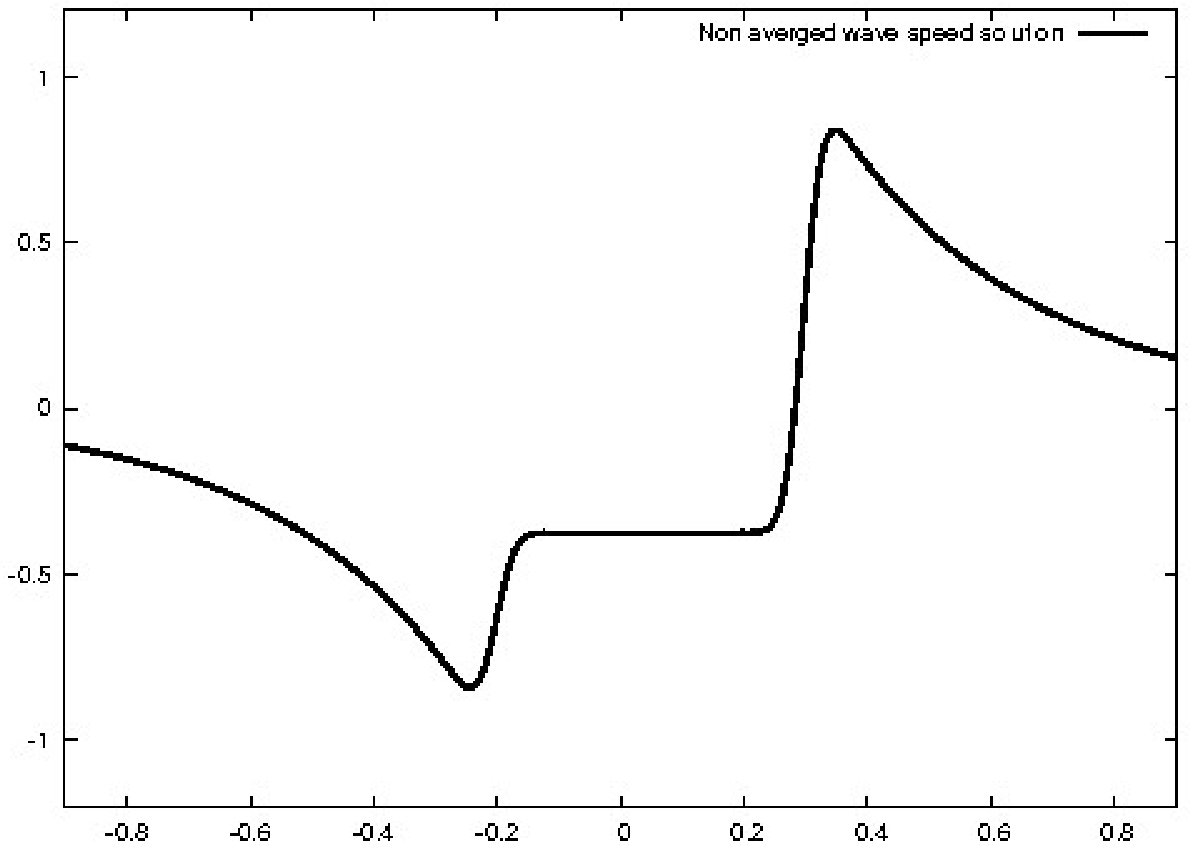} \\
    (a) & & (b) \\
\end{array}$$
\vskip-0.4cm
\caption{(a) Propagation of the initial polynomial condition using the averaged \textit{Riemann} solver, ~(b) Propagation of the initial polynomial condition using the proposed \textit{Riemann} solver}
\label{fig16}
\end{figure}

\begin{figure}[H]
$$\begin{array}{ccc}
    \includegraphics[scale=0.55]{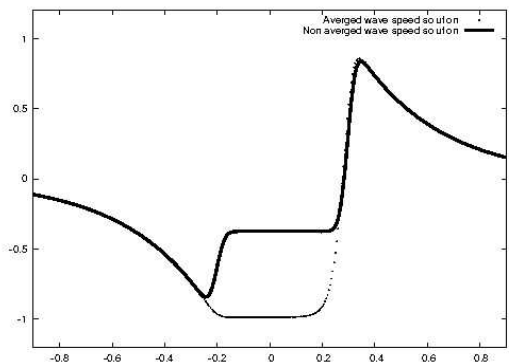} \\
\end{array}$$
\vskip-0.4cm
\caption{Curves of Figure (\ref{fig16}) superimposed}
\label{fig17}
\end{figure}

Now to analyse the case of the discontinuity of the coefficient $a(.)$ and the initial condition at the same location that leads to a product of distribution, we consider the same function $a(.)$ as before, and we add and subtracted a value to the sinusoidal function to create discontinuity at the origin as showed in Figure (\ref{fig12}). The result of the propagation is depicted in Figures (\ref{fig13}, \ref{fig14}), which show a sensitive difference between the two solutions.  Figures (\ref{fig15}, \ref{fig16},  \ref{fig17}) confirm the result for a discontinuous initial polynomial function. as we can see the scheme obtained by waves speed averaging is capable to predict the plateau (the intermediate state) but with a completely wrong value. In another hand this agreed with the theory of generalized function algebra that tries to handle the situation of product of distributions, while the averaging process in the numerical scheme is avoiding this situation.

\section{Conclusions}
The paper dealt with the impact of waves speed averaging process commonly used in \textit{Godunov} type numerical schemes. First a \textit{Riemann} solver for scalar hyperbolic equations with discontinuous coefficients, published first in a conference paper, is presented. After recalling the numerical argument to demonstrate the validity of the solver provided in the conference paper, a new argument based on a regularization procedure and asymptomatic analysis is introduced to reinforce the validity of the proposed solution in absence of a rigorous mathematical proof. The mathematical analysis of the PDE model is briefly discussed pointing out the problem of product of distributions that could results from the discontinuity of the coefficient and the solution at the same location that require the use of \textit{Colombeau} algebra for analyses. To demonstrate the impact of waves speed averaging, a \textit{Gudunov} scheme using the proposed \textit{Riemann} solver is derived  and results are compared to the same scheme using waves speed averaging. The tests show clearly a significant discrepancy of the solutions in the case of discontinuous coefficient with discontinuous solution at the same location. This shows first that the averaging process could lead to a wrong solution and second it makes a connection to the phenomenon of product of distributions. 
Note that this work is a first step toward a construction of a \textit{Riemann} solver for non-linear hyperbolic equations and systems with discontinuous coefficients and its application for instance for inviscid fluxes estimation in the Navier-Stokes and electromagnetic equations discretization.

\section*{Acknowledgement}

This research is supported by the Basque Government 
through the BERC 2014-2017 program and by the 
Spanish Ministry of Economy and Competitiveness 
MINECO: BCAM Severo Ochoa accreditation 
SEV-2013-0323.
Lakhdar Remaki was partially funded by the Project of the Spanish Ministry of Economy and Competitiveness with reference MTM2013-40824-P.

\end{document}